\documentclass[a4paper]{article}
\usepackage[dvips]{graphicx}
\usepackage{amssymb}
\usepackage{amsmath}

\newcommand{\bs}{\boldsymbol}
\makeatletter
\def\dddots{\mathinner{\mkern1mu\raise\p@  
    \hbox{.}\mkern2mu\raise4\p@\hbox{.}\mkern2mu
    \raise7\p@\vbox{\kern7\p@\hbox{.}}\mkern1mu}}%
\makeatother
\newtheorem{theorem}{Theorem}

\newtheorem{Proposition}{Proposition}

\newenvironment{AMS}{\small\bf 2000 AMS subject classification: }{}

\newfont{\BBB}{msbm10 scaled \magstep1}
\newfont{\BBS}{msbm10}
\newcommand{\bb}[1]{\mbox{\BBB #1}}

\begin{document}
\title{\bf Geometric Weakly Admissible Meshes, Discrete Least Squares 
Approximations and Approximate Fekete 
Points\thanks{
Supported by the ``ex-$60\%$'' funds of the University of Padova, 
by the INdAM GNCS and NSERC, Canada.}}

\author{ L. Bos \\ Department of Mathematics and Statistics \\
University of Calgary, Calgary, Alberta \\
Canada T2N 1N4 \\ \\
J.-P. Calvi\\
Institut de Math\'ematiques de Toulouse \\
Universit\'e Paul Sabatier\\
32062, Toulouse Cedex 9\\
France\\
\\
N. Levenberg\\
Department of Mathematics\\
Indiana University\\
Bloomington, Indiana
47405 USA\\
\\
A. Sommariva, M. Vianello\\
Department of Pure and Applied Mathematics\\
University of Padova\\
35121 Padova, Italy }

\maketitle

\begin{abstract}
Using the concept of Geometric Weakly Admissible Meshes (see \S2 below) together with an algorithm based on the classical 
QR factorization of matrices, we compute efficient points for discrete multivariate least squares approximation and Lagrange interpolation.
\end{abstract}

\vskip0.5cm
\noindent
\begin{AMS}
{\rm 41A10, 41A63, 65D05, 65D15, 65Y20.} 
\end{AMS}
\vskip0.2cm
\noindent
{\small{\bf Keywords:} Admissible Meshes, Discrete Least 
Squares Approximation, Approximate Fekete Points, 
Multivariate Polynomial Interpolation.}

\section{Introduction.} 
In a recent paper \cite{CL08}, Calvi and Levenberg have developed the 
theory of ``admissible meshes'' for uniform polynomial approximation over 
multidimensional compact sets. Their theory is for  $\bb{C}^d$, but 
here we restrict our attention to $\bb{R}^d$, the relevant definitions 
and results being easily adaptable to the general case. We begin with 
introducing some notation.
 
Suppose that $K\subset \bb{R}^d$ is compact. We will require that $K$ not 
be too small, that is, that it is
{\it polynomial determining}, i.e., if a polynomial $p$ is zero on $K,$ then it is identically zero. This will
certainly be the case if $K$ contains some open ball, as will be the case for all our examples.

We will let $\bb{P}_n^d$ denote the space of  polynomials of degree at most $n$ in $d$ real variables.

Then, a {\em 
Weakly Admissible 
Mesh} (WAM) of 
$K$   
is a sequence of discrete subsets $A_n\subset K$ such that 
\begin{equation} \label{WAM}
\|p\|_K\leq C(A_n) \|p\|_{A_n}\;,\;\;\forall p\in \bb{P}_n^d\;,
\end{equation}
where ${\rm card}(A_n)$ and the constants $C(A_n)$ have (at most) polynomial 
growth in $n,$ i.e.,
\begin{equation} \label{WAM2}
N:=dim(\bb{P}_n^d)={n+d\choose d}\leq 
{\rm card}(A_n)=\mathcal{O}(n^\alpha)\,\,\alpha>0
\end{equation} 
and 
\begin{equation} \label{WAM3}
C(A_n)=\mathcal{O}(n^\beta)\;,\;\;\beta>0.
\end{equation}
Throughout the paper,  
$\|p\|_X=\max_{\bs{x}\in X}{|p(\bs{x})|}$. 
When $\{C(A_n)\}$ is bounded, $C(A_n)\leq C$, we speak of an {\em 
Admissible 
Mesh} (AM). It is easy to see that (W)AMs satisfy the following 
properties (cf. \cite{CL08}): 
\begin{itemize}

\item for affine mappings of $K$ the image of a WAM is also a WAM, with the same constant $C(A_n)$

\item any sequence of sets containing $A_n$, with polynomially growing 
cardinalities, is a WAM with the same constants $C(A_n)$  

\item any sequence of unisolvent interpolation sets whose Lebesgue 
constant 
grows at most polynomially with $n$ is a WAM, $C(A_n)$ being the Lebesgue 
constant itself

\item a finite union of WAMs is a WAM for the corresponding  
union of compacts, $C(A_n)$ being the maximum of the corresponding 
constants

\item a finite cartesian product of WAMs is a WAM for the corresponding 
product of compacts, $C(A_n)$ being the product of the corresponding 
constants.

\end{itemize}
  
As shown in \cite{CL08}, such meshes are very useful for polynomial 
approximation in 
the max-norm on $K$. In fact, given the classical 
{\em least-square} polynomial approximant on $A_n$ to a function 
$f\in C(K)$, say $\mathcal{L}_n f\in \bb{P}_n^d$, we have that  
\begin{equation} \label{LS} 
\|f-\mathcal{L}_n f\|_K \leq 
\left(1+C(A_n)\left(1+\sqrt{{\rm 
card}(A_n)}\right)\right)\,\min{\{\|f-p\|_K,\,p\in
\bb{P}_n^d\}}\;, 
\end{equation}
see \cite[Thm. 1]{CL08}. 
Moreover, {\em Fekete points} (points that maximize the 
Vandermonde determinant) extracted from a WAM have a 
Lebesgue constant with the bound
\begin{equation} \label{leb_fek}
\Lambda_n\leq N\,C(A_n)\;,
\end{equation} 
that is $C(A_n)$ times the classical bound for the 
Lebesgue constant of true (continuum) Fekete points; see \cite[\S 
4.4]{CL08}. 
Recently, a new algorithm has been proposed
for the computation of Approximate Fekete Points, using 
only standard
tools of numerical linear algebra such as the QR factorization of 
Vandermonde matrices, cf. \cite{BL08,SV08}.    

These results show that in computational applications it 
is important to construct WAMs with 
low cardinalities and slowly increasing constants $C(A_n)$. 
We recall that is always possible to construct easily 
an AM on compact sets that admit a Markov polynomial inequality, 
as is shown in the key result \cite[Thm. 5]{CL08}. 
There are wide classes of such compacts, for example convex bodies, 
and more generally sets satisfying an interior cone condition, 
but the best known bounds for the cardinality of the resulting meshes grow like
$\mathcal{O}(n^{rd})$, where $r$ is the 
exponent of the Markov inequality (typically $r=2$ in the real case).  
This means that, for example, for a 3-dimensional cube 
or ball one should work with $\mathcal{O}(n^6)$ points, 
a number that becomes practically intractable already for 
relatively small values of $n$ (recall that the number of points 
determines the number of rows of the Least Squares (non-sparse) matrices).   
But already in dimension two, working with $\mathcal{O}(n^4)$ 
points makes the construction of polynomial approximants 
at moderate values of $n$  computationally rather expensive. 

The properties of WAMs listed above 
(concerning finite unions and products) suggest an alternative: 
we can obtain 
good meshes, even on complicated geometries, if we are able to 
compute WAMs of low cardinality on standard compact ``pieces''. 
For example, it is immediate to get a  $\mathcal{O}(n^d)$ WAM with 
$C(A_n)=\mathcal{O}(\log^d{(n)})$ 
for the $d$-dimensional cube, 
as the  tensor-product of 1-dimensional Fekete (or Chebyshev) 
points.  

In this paper we introduce the notion of {\em geometric} WAM, 
that is a WAM obtained by a geometric transformation of a suitable 
low-cardinality discretization mesh on some reference compact set, 
like the $d$-dimensional cube. Such geometric WAMs, as well as those 
obtained by finite unions and products, can be used 
directly for discrete Least Squares approximation, as well as for the extraction 
of good interpolation points by means of the Approximate Fekete Points 
algorithm
described in \cite{BL08,SV08}, on compact sets with various geometries. 

In Section 2 we illustrate the idea of geometric WAMs by several 
examples in $\bb{R}^2$,  the disk, triangles, trapezoids, and 
polygons. 
In Section 3 we prove a general result on the asymptotics of approximate 
Fekete points extracted from WAMs by the greedy algorithm in 
\cite{BL08,SV08}.  
Finally, in Section 4 we present some numerical results concerning 
discrete Least Squares approximation and interpolation at Approximate 
Fekete Points on geometric WAMs.

\section{Geometric WAMs.}
Let $K$ and $Q$ be compact subsets of $\bb{R}^d$, 
\begin{equation} \label{t}
\bs{t}:Q\to K\;\;\; \mbox{a {\em surjective} map}
\end{equation}
and $\{S_n\}$ a 
sequence of 
discrete subsets of $Q$. A {\em geometric} WAM of $K$ is a sequence 
\begin{equation} \label{geowam}
A_n=\bs{t}(S_n)\;, 
\end{equation}
where the defining properties of a WAM,  (\ref{WAM})-(\ref{WAM3}) stem from 
the  
``geometric structure'' of  $K$, $Q$, $\bs{t}$, and $S_n$.
To make more precise this still somewhat vague notion, we give some 
illustrative examples in $\bb{R}^2$, where the reference 
compact $Q$ is a rectangle.

\subsection{The disk.}
A geometric WAM of the disk can be immediately obtained by 
working with {\em polar coordinates}, i.e., by considering the map
\[
\bs{t}: Q=[0,1]\times [0,2\pi] \to K=\{\bs{x}:\,\|\bs{x}\|_2\leq 1\}
\]
\begin{equation} \label{polar}
\bs{t}(r,\phi)=(r\cos{\phi},r\sin{\phi})\;.
\end{equation}
We now state and prove the following
\begin{Proposition}
The sequence of polar grids 
\[
S_n=\{(r_j,\phi_k)\}_{j,k}=\left\{\frac{1}{2}+\frac{1}{2}\,
\cos{\frac{j\pi}{n}}\,,\,0\leq 
j\leq 
n\right\} \times 
\left\{\frac{2 \pi 
k}{2n+1}\,,\,0\leq k\leq 2n\right\}
\]
gives a WAM $A_n=\bs{t}(S_n)$ of the unit disk,  
such that $C(A_n)=\mathcal{O}(\log^2{(n)})$ and ${\rm card}(A_n)=2n^2+n+1$.

\end{Proposition}
\vskip0.2cm
\noindent
{\em Proof.} 
Observe that given a polynomial $p\in \bb{P}_n^2$, 
when restricted to the disk in polar coordinates, 
$q(r,\phi)=p(\bs{t}(r,\phi)),$ becomes a polynomial of degree $n$  
in $r$ for any fixed $\phi$, and a trigonometric polynomial of degree $n$ 
in $\phi$ for any fixed $r$. Recall that $n+1$ Chebyshev-Lobatto points 
are near-optimal for 1-dimensional polynomial interpolation, and 
$2n+1$ equally 
spaced points are near-optimal for trigonometric interpolation, both 
having a  
Lebesgue constant $\mathcal{O}(\log{(n)})$; cf. \cite{DDP81,R69}.   
Now, for every $p\in \bb{P}_n^2$ we can write
\[
|p(x_1,x_2)|= |q(r,\phi)|=|p(r \cos{\phi},r \sin{\phi}|
\leq c_1\log{(n)}\,\max_j{|q(r_j,\phi)|}
\]
where $c_1$ is independent of $\phi.$ 
since the $\{r_j\}$ are the $n+1$ Chebyshev-Lobatto points in $[0,1].$ 
Further
\[
|q(r_j,\phi)| \leq c_2 \log{(n)}\,\max_k{|q(r_j,\phi_k)|} 
\]
where $c_2$ is independent of $j,$
since the $\{\phi_k\}$ are $2n+1$ equally spaced points in 
$[0,2\pi].$  Thus 
\[
|p(x_1,x_2)|\leq c_1 c_2 
\log^2{(n)}\,\max_{j,k}{|q(r_j,\phi_k)|}\;,\;\;\forall (x_1,x_2)\in K\;,
\]
i.e., $A_n=\bs{t}(S_n)=\{(r_j \cos{\phi_k},r_j \sin{\phi_k})\}$ is a WAM for
the disk
with $C(A_n)=O(\log^2{(n)})$. We conclude by observing that the number of 
distinct points in $\bs{t}(S_n)$ is, due to the fact that $r_n=0,$
${\rm card}(S_n)-2n=(n+1)(2n+1)-2n=2n^2+n+1$.$\;\;\;\square$
\vskip0.2cm
\noindent

\begin{figure}[!ht]
\centering
\includegraphics[scale=0.38,clip]{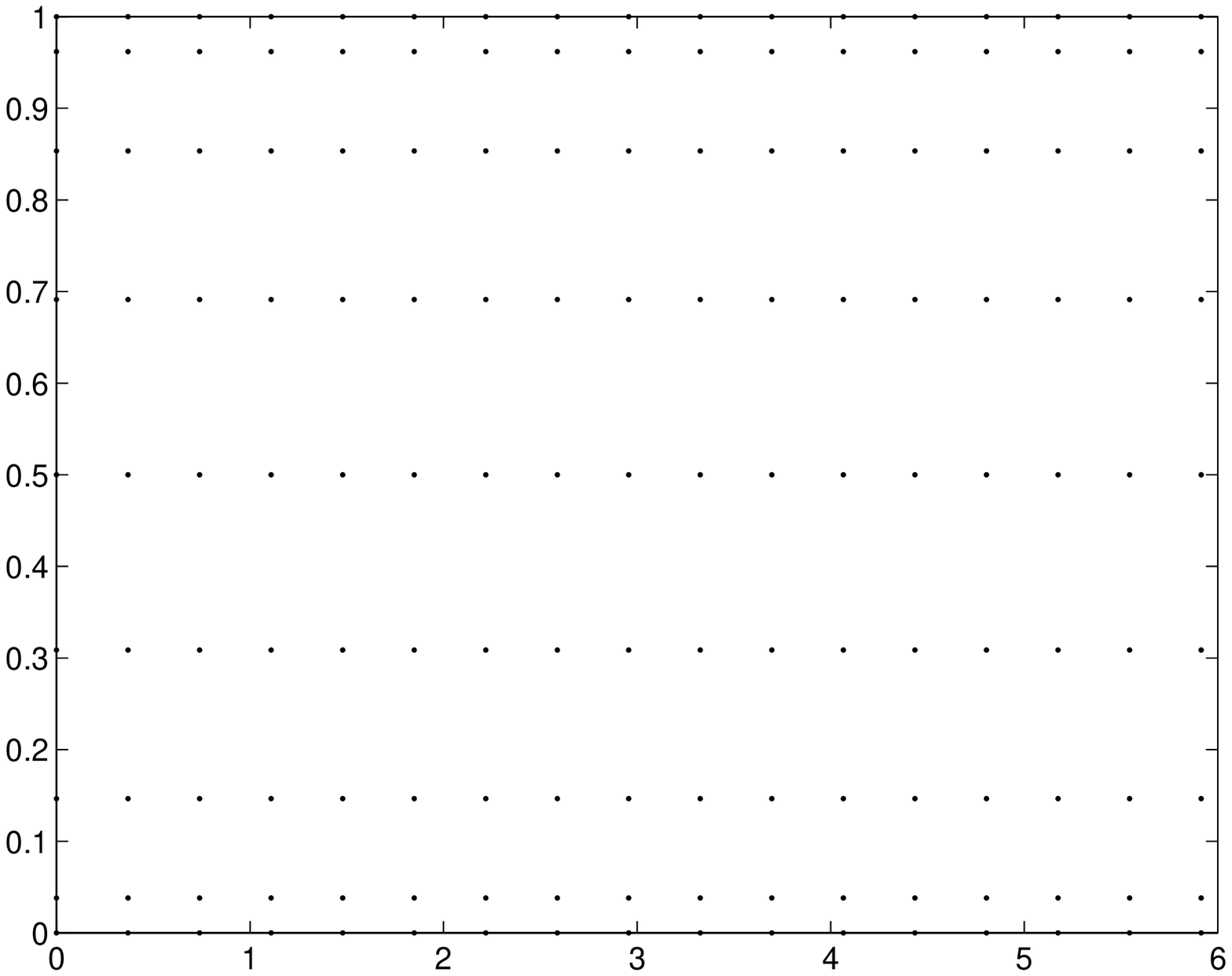}\hfill
\includegraphics[scale=0.38,clip]{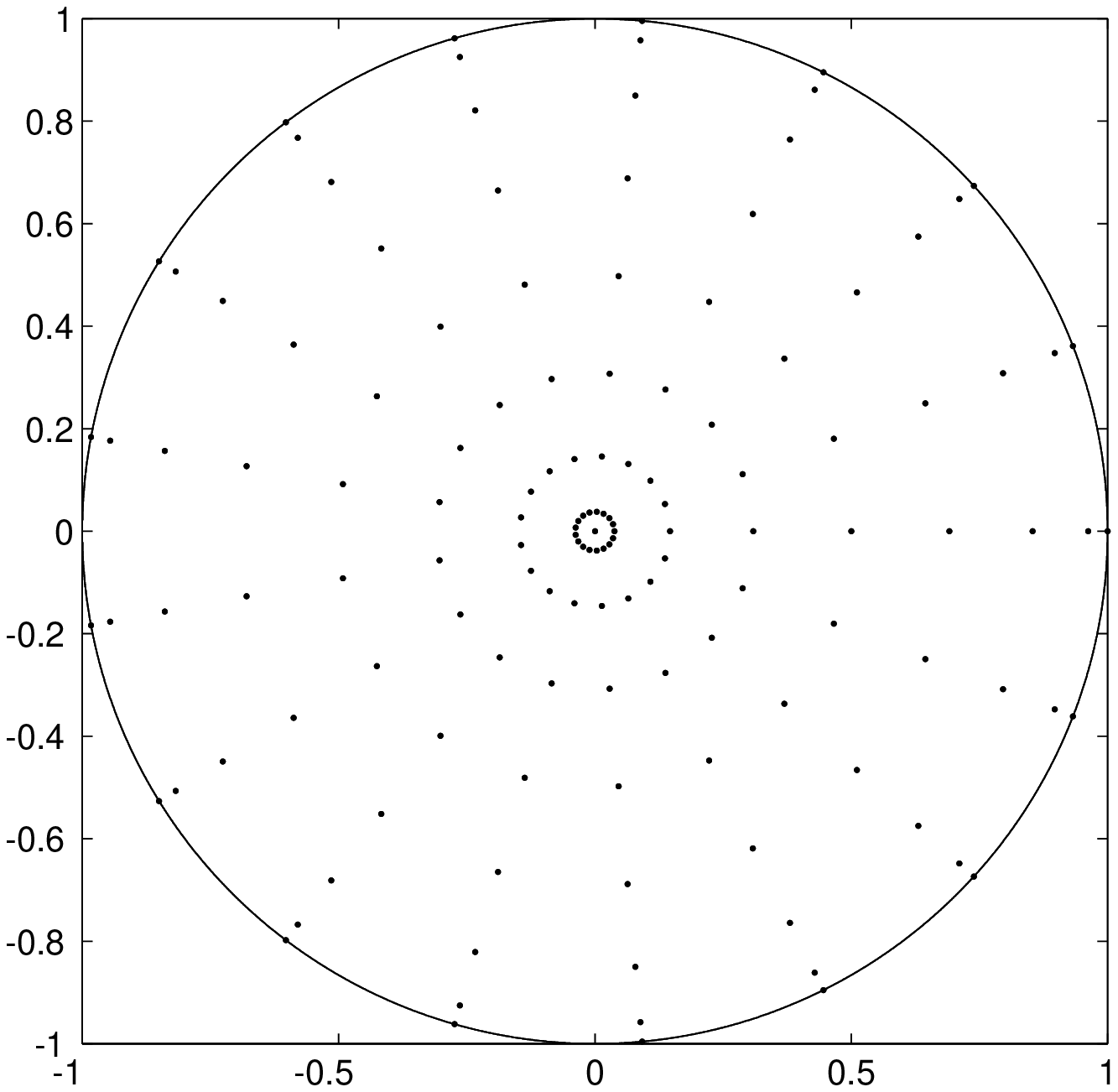}\hfill
\caption{The polar grid and the corresponding geometric 
WAM of degree $n=8$ for the unit disk.}
\label{disk}
\end{figure}

In Figure \ref{disk}, we display the polar grid $S_n$ and the WAM $A_n$
in the unit disk for $n=8$
($153$ points and $137$ points, respectively). In view of the structure 
of $S_n$ and $\bs{t}$, the points of the geometric WAM cluster at the 
boundary and at the center of the disk. Note that this technique can also
be used to construct geometric WAMs for annuli and even a ball in higher dimensions.
However we will not pursue this here.

\subsection{Polynomial maps.}
In the case that $\bs{t}$ is a {\em polynomial map}, we can give the 
following general 

\begin{Proposition}
Let $\{B_n\}$ be a WAM of $Q$, and $\bs{t}:Q\to K$ a surjective 
polynomial map 
of 
degree $k$, i.e., 
$\bs{t}(\bs{y})=(t_1(\bs{y}),...,t_d(\bs{y}))$ 
with $t_j\in \bb{P}^d_k$. 
Then $A_n=\bs{t}(S_n)$, with $S_n=B_{kn}$, is a WAM of $K$ such that  
$C(A_n)=C(B_{kn})$.  
\end{Proposition}
\vskip0.2cm
\noindent
{\em Proof.} Observing that for every  
$\bs{x}\in K$ there exists some $\bs{y}\in Q$ such that 
$\bs{x}=\bs{t}(\bs{y})$, 
and that  
for every $p\in \bb{P}_n^d$ the composition $p(\bs{t}(\cdot))$ 
is a polynomial of degree $\leq kn$, we have   

$$
|p(\bs{x})|=|p(\bs{t}(\bs{y}))|\leq 
C(B_{kn})\|p(\bs{t}(\cdot))\|_{B_{kn}}= 
C(B_{kn})\|p\|_{\bs{t}(B_{kn})}\;.\;\;\;\square
$$

\subsubsection{Triangles.}
We begin by observing that for the square $Q=[-1,1]^2$ several WAMs 
$B_n$ formed by good interpolation points are known, namely
\begin{itemize}

\item tensor-products of 1-dimensional near-optimal interpolation 
points, such as Chebyshev-Lobatto points, Gauss-Lobatto points, and 
others; 

\item the {\em Padua points}, recently studied in 
\cite{BCDMVX06,BDMVX07,CDMV05,CDMV08a,CDMV08b}, that are near-optimal 
for total-degree interpolation.

\end{itemize}

All these WAMs have $C(B_n)=\mathcal{O}(\log^2{(n)})$, but 
the Padua points have minimal cardinality 
$N={\rm dim}(\bb{P}^2_n)=(n+1)(n+2)/2$, 
whereas the cardinality of 
tensor-product points is $(n+1)^2={\rm dim}(\bb{P}^1_n \bigotimes   
\bb{P}^1_n)$. We recall that, given the 1-dimensional 
Chebyshev-Lobatto points 
\begin{equation} \label{cheblob}
\mathcal{C}_{n+1}=\{\cos{(j\pi/n)},\,0\leq j \leq n\}\;, 
\end{equation} 
the Padua points of degree $n$ are given by the union of two 
Chebyshev-Lobatto grids
\begin{equation} \label{padua}
B_n=\mbox{Pad}_n=(\mathcal{C}^{\mbox{\footnotesize{odd}}}_{n+1} \times 
\mathcal{C}^{\mbox{\footnotesize{even}}}_{n+2}) \cup 
(\mathcal{C}^{\mbox{\footnotesize{even}}}_{n+1} 
\times
\mathcal{C}^{\mbox{\footnotesize{odd}}}_{n+2}) \subset \mathcal{C}_{n+1} 
\times 
\mathcal{C}_{n+2}\;. 
\end{equation}

\begin{figure}[!ht]
\centering
\includegraphics[scale=0.42,clip]{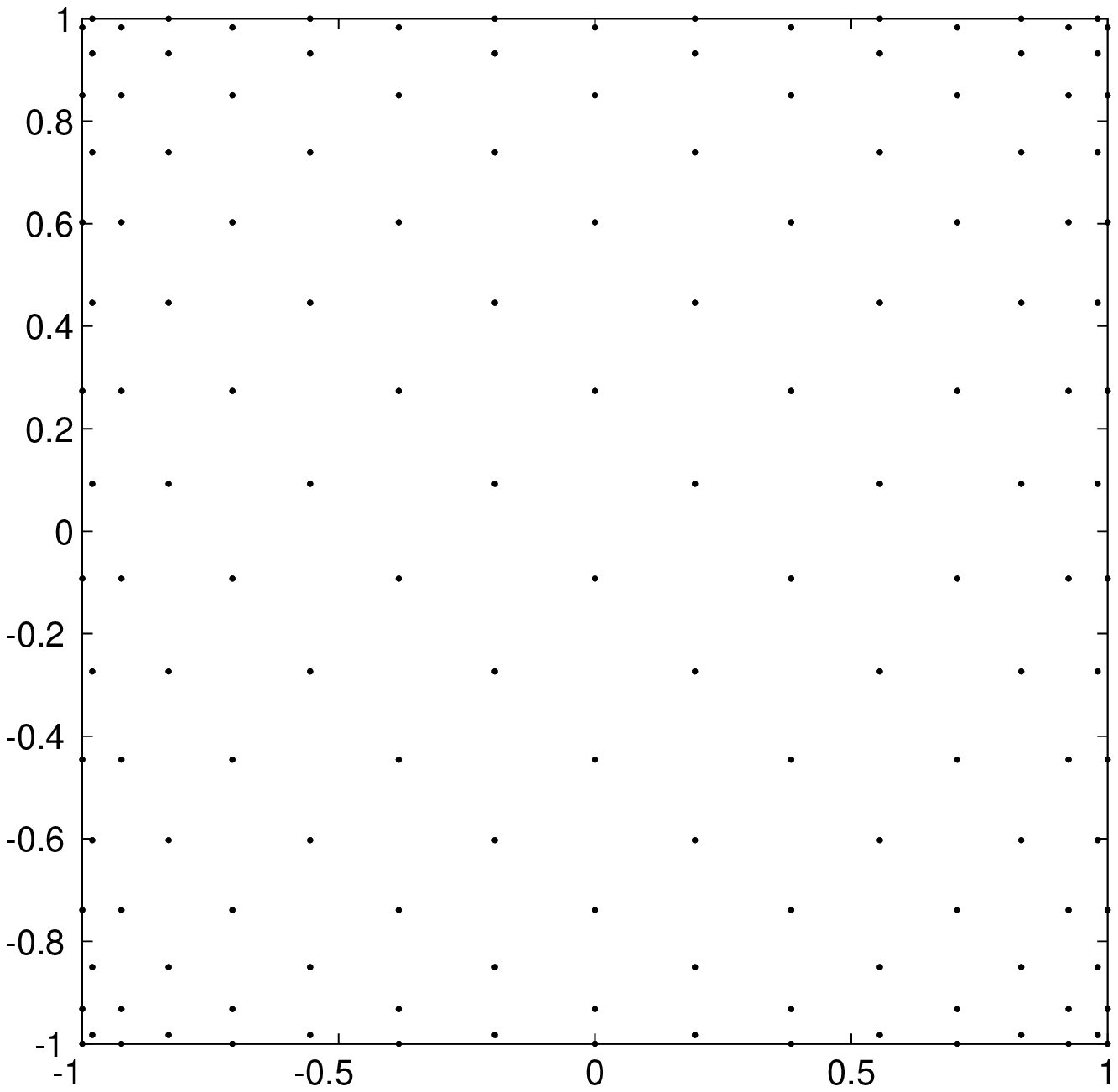}\hfill
\includegraphics[scale=0.42,clip]{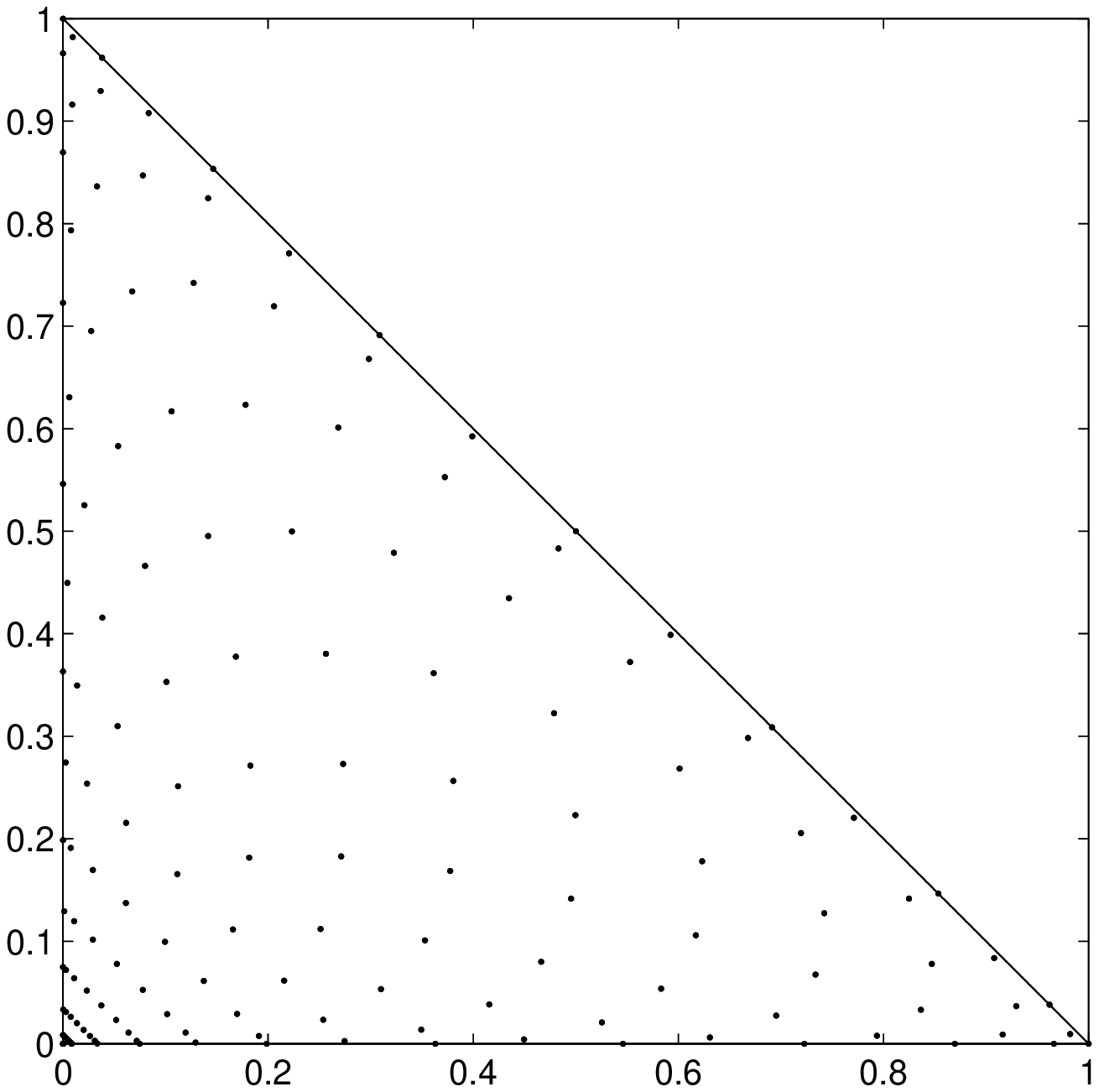}\hfill
\caption{The Padua points of degree $2n=16$ and the corresponding 
geometric
WAM of degree $n=8$ for the unit simplex.}
\label{simplex}
\end{figure}

There is a simple quadratic map from the square to any triangle 
with vertices $\bs{u}=(u_1,u_2)$,  $\bs{v}=(v_1,v_2)$,  
$\bs{w}=(w_1,w_2)$,  
namely the Duffy transformation (cf. \cite{D82}) 
\begin{equation} \label{duffy}
\bs{t}(\bs{y})=\frac{1}{4}\,(\bs{v}-\bs{u})(1+y_1)(1-y_2)
+\frac{1}{2}\,(\bs{w}-\bs{u})(1+y_2)+\bs{u}\;,
\end{equation}
which collapses one side of the square (here $y_2=1$) onto 
a vertex of the triangle (here $\bs{w}$). By this map, the Padua points 
$B_{2n}$ of the square (cf. (\ref{padua})) are transformed to the 
WAM $A_n=\bs{t}(B_{2n})$ of the triangle, with constants 
$C(A_n)=\mathcal{O}(\log^2{(2n)})=\mathcal{O}(\log^2{(n)})$. 
The number of distinct points in $\bs{t}(B_{2n})$ is 
${\rm card}(A_n)={\rm card}(B_{2n})
-{\rm card}(\mathcal{C}^{\mbox{\footnotesize{odd}}}_{2n+1})+1
={\rm dim}(\bb{P}_{2n}^2)-(n+1)
=2n^2+2n$. 
In Figure \ref{simplex}, we display the WAMs $B_{2n}$ in the square and 
$A_n$ 
in the unit simplex for $n=8$ 
($153$ points and $144$ points, respectively). In view of the properties 
of the Padua points, the points of the geometric WAM cluster at the sides and 
especially at the vertices of the simplex.   

\subsubsection{Polynomial trapezoids.}
We consider here bidimensional compact sets of the form
\begin{equation} \label{trap}
K=\{\bs{x}=(x_1,x_2):\,a\leq x_1\leq b,\,g_1(x_1)\leq x_2\leq 
g_2(x_1)\}\;, 
\end{equation}
where $g_1,g_2\in \bb{P}^1_\nu$. 
A polynomial map $\bs{t}$ of degree $k=\nu+1$ from the  square 
$Q=[-1,1]^2$ onto $K$ is 
\begin{equation} \label{trapmap}
t_1(\bs{y})=\frac{b-a}{2}y_1+\frac{b+a}{2}\;,\;\;
t_2(\bs{y})=\frac{g_2(t_1)-g_1(t_1)}{2}y_2+\frac{g_2(t_1)+g_1(t_1)}{2}\;.
\end{equation}
Observe that a triangle could be treated (up to an affine tranformation) 
as a degenerate linear trapezoid. 

By Proposition 2, the Padua points $B_{kn}$ of the square are  
mapped to the WAM $A_n=\bs{t}(B_{kn})$ of the polynomial trapezoid, 
with $C(A_n)=\mathcal{O}(log^2{(kn)})$ and ${\rm card}(A_n)\leq {\rm card}(B_{kn})
=dim(\bb{P}_{kn}^2)=(kn+1)(kn+2)/2$.  

In Figure \ref{trap} we show the WAMs $A_n=\bs{t}(B_{2n})$ and 
$A_n=\bs{t}(B_{4n})$ 
obtained by mapping the Padua points onto a linear trapezoid (quadratic 
map) and a cubic trapezoid (quartic map), again for $n=8$ 
(the numbers of points are $231={\rm dim}(\bb{P}_{20}^2)$ and $561={\rm 
dim}(\bb{P}_{32}^2)$,  
respectively). As expected, we observe clustering at the sides and at the 
vertices of the trapezoids. Notice that it could happen that 
${\rm card}(A_n)<{\rm card}(B_{kn})$, 
namely when the graphs of $g_1$ and $g_2$ intersect
at some $x_1\in t_1(\mathcal{C}_{kn+1})$ (the $kn+1$ Chebyshev-Lobatto 
points of $[a,b]$).

\begin{figure}[!ht]
\centering
\includegraphics[scale=0.42,clip]{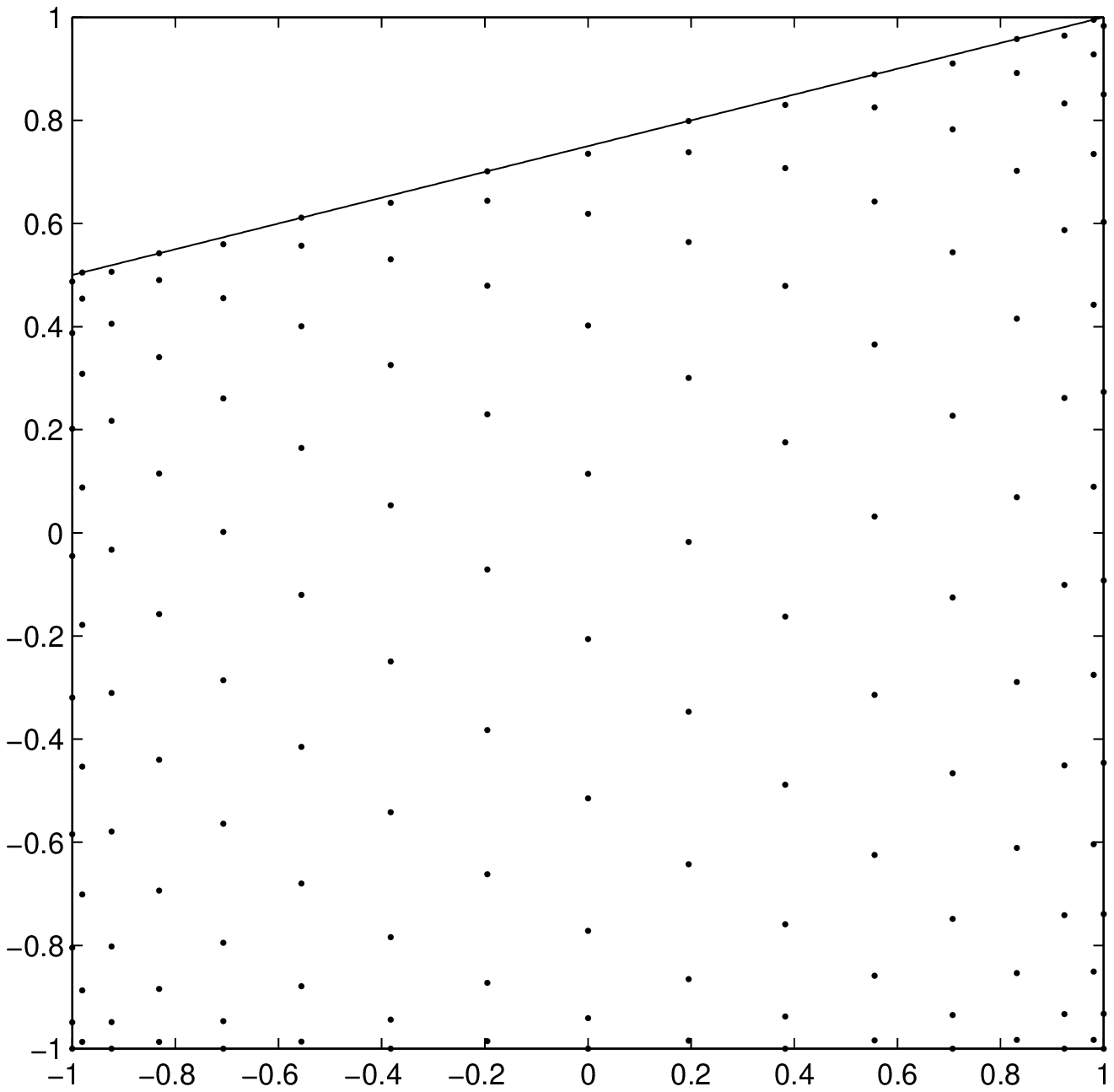}\hfill
\includegraphics[scale=0.42,clip]{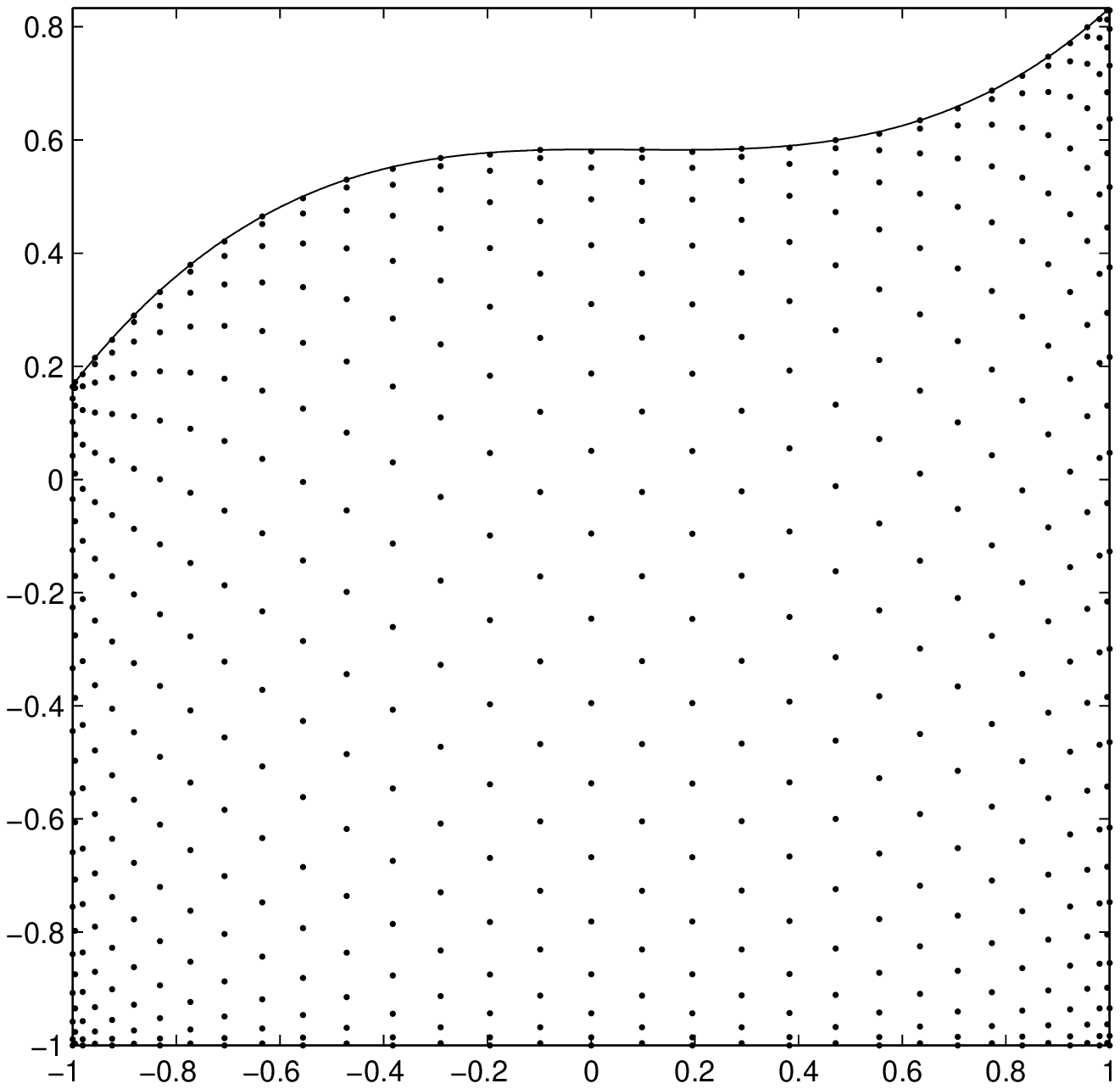}\hfill
\caption{Examples of geometric WAMs of degree $n=8$ 
for a linear and a cubic trapezoid.}
\label{trap3}
\end{figure}

\subsubsection{Finite unions: polygons.}
A relevant example concerning finite unions is given by polygons, 
which are very important in applications such as
2-dimensional computational geometry. As known, any simple 
(no interlaced sides) 
and simply connected polygon with $m$ vertices can be subdivided into 
$m-2$ triangles, and this can be done by fast algorithms, cf. e.g. 
\cite{H01,NM95}. Once this rough triangulation is at hand, we can 
immediately obtain a geometric WAM $A_n$ for the polygon by union of the 
geometric 
WAMs constructed by mapping the Padua points on the triangles, with 
$C(A_n)=\mathcal{O}(\log^2{(n)})$ 
and ${\rm card}(A_n)\leq (m-2)(2n^2+2n)$, in view of the basic properties of 
WAMs and the results of Section 2.2.1. The points of the union WAM 
will cluster especially at the triangles common sides and vertices. 

\begin{figure}[!ht]
\centering
\includegraphics[scale=0.42,clip]{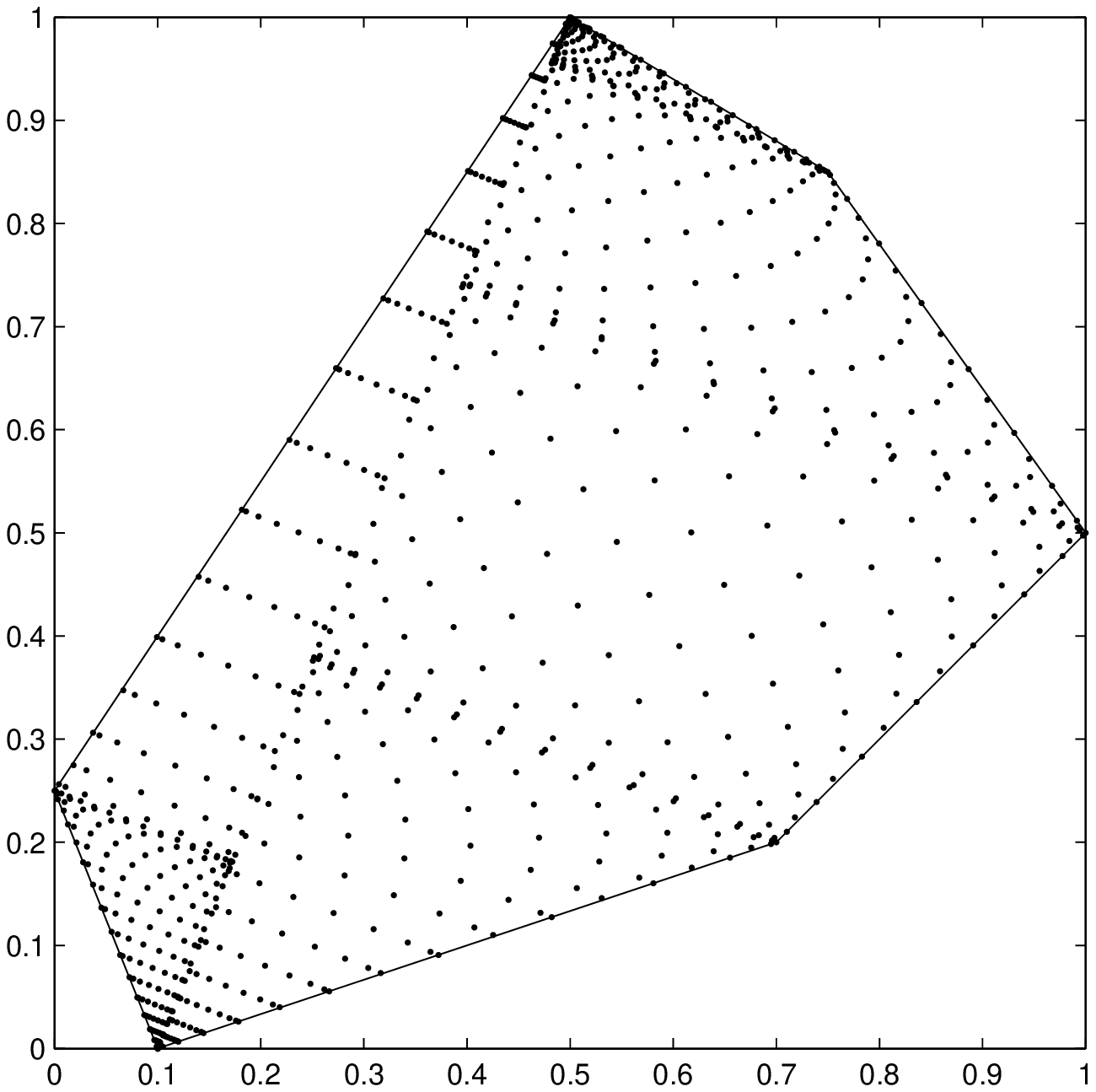}\hfill   
\includegraphics[scale=0.42,clip]{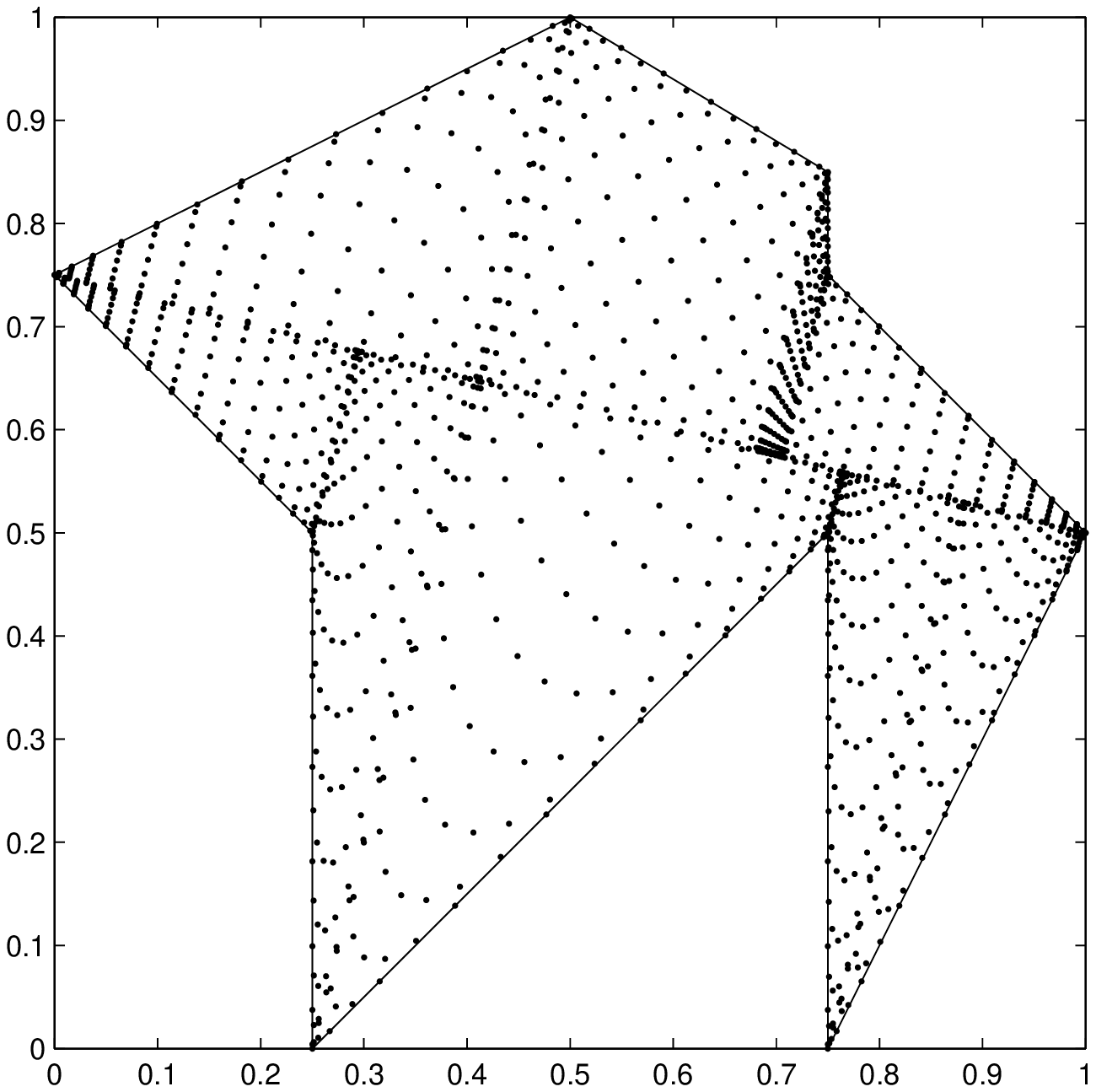}\hfill   
\caption{Examples of geometric WAMs of degree $n=8$
for a convex and a nonconvex polygon.}
\label{polygons}
\end{figure}

A similar approach is to subdivide into 
linear trapezoids, where we again obtain by finite union a geometric WAM $A_n$ 
for the polygon with $C(A_n)=\mathcal{O}(\log^2{(n)})$
and ${\rm card}(A_n)=\mathcal{O}(mn^2)$ (see Section 2.2.2). 
In Figure \ref{polygons} 
we show geometric WAMs generated by subdivision into linear trapezoids 
of a convex and a nonconvex polygon, for $n=8$. The method adopted, 
which works for a wide class of polygons, is that 
used for the generation of algebraic cubature points in \cite{SV07},  
where the trapezoidal panels are obtained simply by orthogonal 
projection of the 
sides on a fixed reference line (observe that the points cluster at the 
sides and at the line).       

\section{Approximate Fekete Points from WAMs.}

In this section we go to the general setting of polynomial interpolation  
in several complex variables. 
Suppose that $K\subset\bb{C}^d$ is a compact polynomially determining set.  
If ${\cal B}_n=\{p_1,p_2,\cdots,p_N\}$ is a basis for $\bb{P}_n^d,$ where
$N:={\rm dim}(\bb{P}_n^d)$ and  $Z_n\subset K$ is a discrete subset of 
cardinality $N,$ 
then
\begin{equation}\label{VdmDet}
{\rm vdm}(Z_n;{\cal B}_n):={\rm det}([p(\bs{a})]_{p\in {\cal 
B}_n,\bs{a}\in Z_n})
\end{equation}
is called the corresponding Vandermonde determinant. If the determinant is non-zero then we may
form the so-called fundamental Lagrange interpolating polynomials
\begin{equation}\label{LagrangePoly}
\ell_{\bs{a}}(\bs{z}):={{\rm vdm}(Z_n\backslash\{\bs{a}\}\cup 
\{\bs{z}\};{\cal B}_n)\over 
{\rm vdm}(Z_n;{\cal B}_n)}. 
\end{equation}
Indeed, then for any $f\,:\,K\to\bb{C}$ the polynomial (of degree $n$)
\begin{equation}\label{LagrangeForm}
\Pi_n(\bs{z})=\sum_{\bs{a}\in Z_n}f(\bs{a})\ell_{\bs{a}}(\bs{z})
\end{equation}
interpolates $f$ at the points of $Z_n,$ i.e., $\Pi_n(\bs{a})=f(\bs{a}),$ 
$\bs{a}\in Z_n.$
A set of  points $F_n\subset K$ which maximize ${\rm vdm}(Z_n;{\cal B}_n)$ 
as a function of $Z_n,$ are called Fekete points of degree $n$ for $K$ and have the
special property that $\|\ell_{\bs{a}}\|_K=1$ 
(and hence Lebesgue constant bounded by $N$) and provide a very good 
(often excellent) set of interpolation points
for $K.$ However, they are typically very difficult to compute, even 
for moderate values of $n.$

For any $\bb{P}_n^d$ determining subset $A_n\subset K$ (thought of as 
a sufficiently good discrete model of $K$) the algorithm introduced 
by Sommariva and Vianello in \cite{SV08} and studied by Bos and 
Levenberg in \cite{BL08} selects, in a simple and 
efficient manner, a subset ${\cal F}_n\subset A_n$ of Approximate Fekete 
Points, 
and hence provides
a practical alternative to {\it true} Fekete points, $F_n.$
The optimization problem is nonlinear, and large-scale already for 
moderate values of $n$, but the algorithm
is able to give an approximate solution using only standard tools of
numerical linear algebra.

We sketch here the algorithm, in a Matlab-like
notation. The goal is extracting a maximum volume
square submatrix
from the rectangular $N \times \mbox{card}(A_n)$ Vandermonde matrix
\begin{equation} \label{vander}
V=[p(\bs{a})]_{p\in {\cal B}_n,{\bs{a}}\in A_n}\;, 
\end{equation}
where the polynomial basis and the array of points have been 
(arbitrarily) ordered.
The core is given by the following
iteration:
\vskip0.3cm
\noindent
{\bf Algorithm greedy} (max volume submatrix of a matrix $V\in
\mathbb{R}^{N\times
M}$, $M>N$)
\begin{itemize}
\item $\hspace{0.1cm}ind=[\;]\;;$
\item
$\hspace{0.1cm} \mbox{for}\;k=1,\dots,N$
\begin{itemize}
\item ``select the largest norm column $col_{i_k}(V)$'';
$ind=[ind,i_k]$;
\item ``remove from every column
of $V$ its orthogonal projection onto $col_{i_k}$;
\end{itemize}
$\hspace{0.1cm} \mbox{end}$;
\end{itemize}
\vskip0.2cm
\noindent
which works when $V$ is full rank and gives an approximate
solution to the NP-hard maximum volume problem; cf. \cite{CMI07}. 
Then, we can extract the Approximate Fekete Points 
$$
{\cal F}_n=A_n(ind)=(A_n(i_1),\dots,A_n(i_N))\;.
$$
The algorithm can be
conveniently implemented by
the well-known QR
factorization with column pivoting, originally proposed by Businger and
Golub in 1965 \cite{BG65}, and used for example by the Matlab ``mldivide''
or ``$\backslash$'' operator in the solution of underdetermined linear
systems (via the LAPACK routine DGEQP3, cf. \cite{lapack,matlab}).

Some remarks on the polynomial basis ${\cal B}_n$ are in order. First note
that if ${\cal B}'_n:=\{q_1,\cdots,q_N\}$ is some other basis of 
$\bb{P}_n^d$ 
then there is a transition matrix $T_N\in\bb{C}^{N\times N}$ so that
${\cal B}'_n={\cal B}_nT_N$. It is easy to verify that then
\[{\rm vdm}(Z_n;{\cal B}'_n)={\rm det}(T_N){\rm vdm}(Z_n;{\cal B}_n).\]
Hence {\it true} Fekete points do not depend on the basis used and also
the Lagrange polynomials $\ell_{\bs{a}}$ of (\ref{LagrangePoly}) are
independent of the basis. Moreover, if $A_n$ and $B_n$ are two point sets for
which
\[|{\rm vdm}(A_n,{\cal B}_n)|\ge C |{\rm vdm}(B_n,{\cal B}_n)|\]
for some constant $C,$ then the same inequality holds using the basis 
${\cal B}'_n,$
i.e.,
\[ |{\rm vdm}(A_n,{\cal B}'_n)|\ge C |{\rm vdm}(B_n,{\cal B}'_n)|\]
since both sides scale by the same factor.

\medskip
The greedy Algorithm described above {\it is} in general 
affected by the basis. But
this not withstanding,
in the theorem below we show that if the
initial set $A_n\subset K$ is a WAM, 
then the so selected approximate 
Fekete points, using any basis, and the true Fekete points for $K$
both have the same asymptotic distribution. 

\begin{theorem}\label{ConvgThm} Suppose that $K\subset \bb{C}^d$ is 
compact, non-pluripolar, polynomially convex and regular (in the sense of Pluripotential theory) and that for $n=1,2,\cdots,$ $A_n\subset K$ is a WAM. 
Let $\{\bs{b}_1^{(n)},\bs{b}_2^{(n)},\cdots,\bs{b}_N^{(n)}\}$ be the 
Approximate Fekete Points selected from $A_n$ by
the greedy Algorithm of \cite{SV08}, using any basis ${\cal B}_n.$ We 
denote by  $m_n$ 
the sum of the degree of the
$N$ monomials of degree at most $n,$ i.e., $m_n=dnN/(d+1).$ Then

(1) $\displaystyle{\lim_{n\to \infty} 
|{\rm vdm}|(\bs{b}_1^{(n)},...,\bs{b}_N^{(n)})|^{1/m_n}=\tau(K),}$ the 
transfinite diameter of $K$

\noindent and

(2) the discrete probability measures $\mu_n:= 
{1\over N}\sum_{j=1}^{N} \delta_{\bs{b}_j^{(n)}}$ converge weak-* to the 
pluripotential-theoretic {\it equilibrium measure} $d\mu_K$ of $K$.
\end{theorem}

\noindent {\bf Remark.} For $K=[-1,1]$, 
$d\mu_{[-1,1]} = {1\over \pi}{1\over \sqrt {1-x^2}}dx$; 
for $K$ the unit circle $S^1$, $d\mu_{S^1} ={1\over 2\pi}d\theta$. If $K\subset \bb{R}^d\subset\bb{C}^d$ is compact,
then $K$ is automatically polynomially convex.
We refer the reader to \cite{K91} for other examples and more 
on complex pluripotential theory. 

\medskip
By
\[|{\rm vdm}(\bs{z}_1^{(n)},...,\bs{z}_N^{(n)})|\]
we will mean the Vandermonde determinant computed using the standard monomial basis.

Note also that a set of true Fekete points $F_n$ is also a WAM and 
hence we may take $A_n=F_n,$ in which case
the algorithm will select $B_n=F_n$ (there is no other choice) 
and so the true Fekete points must necessarily also
have these two properties.

\noindent {\bf Proof.} We suppose that 
$F_n=\{\bs{f}^{(n)}_1,\bs{f}^{(n)}_2,\cdots,\bs{f}^{(n)}_N\}\subset K$ is 
a
set of true Fekete points for $K.$  
Suppose further that 
$\{\bs{t}^{(n)}_1,\bs{t}^{(n)}_2,\cdots,\bs{t}^{(n)}_N\}\subset 
A_n$ is a set of true Fekete points of degree
$n$ for $A_n$ and that $\ell_i^{(n)}$ are the corresponding Lagrange polynomials. Then,
\[ \|\ell_i^{(n)}\|_K\le C(A_n)\|\ell_i^{(n)}\|_{A_n}=C(A_n).\]
It follows that the associated Lebesgue constants
\[\Lambda_n:=\max_{\bs{z}\in K}\sum_{i=1}^N|\ell_i^{(n)}(\bs{z})|\le 
NC(A_n)\]
and hence, since $C(A_n)$ is of polynomial growth, 
\[\lim_{n\to\infty}\Lambda_n^{1/n}=1.\]
By Theorem 4.1 of \cite{BBCL92} 
\begin{equation}\label{true4An}
\lim_{n\to\infty} |{\rm 
vdm}(\bs{t}_1^{(n)},...,\bs{t}_N^{(n)})|^{1/m_n}=\tau(K).
\end{equation}
By Zaharjuta's famous result \cite{Z}, we also have
\begin{equation}\label{true4K}
\lim_{n\to\infty} |{\rm 
vdm}(\bs{f}_1^{(n)},...,\bs{f}_N^{(n)})|^{1/m_n}=\tau(K).
\end{equation}
Further, by \cite{CMI07}, we have (cf. the remarks on bases preceeding the statement of the Theorem)
\[ |{\rm vdm}(\bs{b}_1^{(n)},...,\bs{b}_N^{(n)})|\ge {1\over N!} |{\rm 
vdm}(\bs{t}_1^{(n)},...,\bs{t}_N^{(n)})|\]
and hence
\begin{eqnarray*}
|{\rm vdm}(\bs{f}^{(n)}_1,\bs{f}^{(n)}_2,\cdots,\bs{f}^{(n)}_{N})|&\ge& 
|{\rm vdm}(\bs{b}^{(n)}_1,\bs{b}^{(n)}_2,\cdots,\bs{b}^{(n)}_{N})| \\
&\ge& {1\over N!}|{\rm 
vdm}(\bs{t}^{(n)}_1,\bs{t}^{(n)}_2,\cdots,\bs{t}^{(n)}_{N})|.
\end{eqnarray*}
Thus by (\ref{true4An}) and (\ref{true4K}) we have
\[\lim_{n\to\infty}|{\rm 
vdm}(\bs{b}_1^{(n)},\bs{b}_2^{(n)},\cdots,\bs{b}_N^{(n)})|^{1/m_n}=\tau(K)\]
as 
\[\lim_{n\to\infty} (N!)^{1/m_n}=1.\]

The final statement, that $\mu_n$ converges weak-* to $d\mu_K$ then  follows by the main result of Berman and Bouksom \cite{BB08} (see also \cite{BBLW08}). $\square$

\section{Numerical results.}
In this section, we present a suite of numerical tests concerning 
discrete least squares approximation on geometric WAMs and polynomial 
interpolation at Approximate Fekete Points extracted from them. The 
tests concern the 2-dimensional compact sets discussed  
in Section 2, that are the unit disk, the unit simplex, a linear and 
a cubic trapezoid, a convex and a nonconvex polygon; see Figures 
\ref{disk}-\ref{polygons}. All the tests have been done in 
Matlab (see \cite{matlab}), by an Intel-Centrino Duo T-2400 processor with 
1 Gb RAM.  

In order to compute the Approximate Fekete Points, we have actually used 
a refined version of the greedy algorithm of the previous section, which 
is sketched below. 
\vskip0.3cm  
\noindent  
{\bf Algorithm greedy with iterative QR refinement of the basis} 
\begin{itemize}
\item take the Vandermonde matrix $V$ in (\ref{vander});
\item $\hspace{0.1cm} V_0=V^t\;;\;T_0=I\;;$
\item $\hspace{0.1cm} \mbox{for}\;\;k=0,\dots,s-1$
\begin{itemize}
\item[$\,$] $\hspace{0.1cm}
V_k=Q_kR_k\;;\;U_k=\mbox{inv}(R_k)\;;$   
\item[$\,$] $\hspace{0.1cm}
V_{k+1}=V_kU_k\;;\;T_{k+1}=T_kU_k\;;$ 
\end{itemize}
$\hspace{0.2cm}$\mbox{end}$\;;$
\item $\;b=(1,\dots,1)^t\;; \hspace{0.2cm}$ (the choice of
$b$ is
irrelevant in practice)
\item $\hspace{0.1cm} w=V_s^t\backslash b\;;\;ind=\mbox{find}
(w\neq 0)\;;\;
{\cal F}_n=A_n(ind)\;;$\\
\end{itemize}
\vskip0.2cm
\noindent
The greedy algorithm is implemented directly by the last row above 
(in Matlab), irrespectively of the actual value of the
vector $b$, and produces a set of Approximate Fekete Points ${\cal F}_n$. 
The
{\em for} loop above implements a change of polynomial basis from
$(p_1,\dots,p_N)$
to the nearly-orthogonal basis
$(q_1,\dots,q_N)=(p_1,\dots,p_N) T_s$ with respect to the 
discrete inner product $(f,g)=\sum_{\bs{a}\in A_n}{f(\bs{a})\,g(\bs{a})}$,
whose main aim is to cope possible numerical rank-deficiency and severe
ill-conditioning arising with nonorthogonal bases (usually
$s=1$ or $s=2$ iterations suffice); for a complete
discussion of this algorithm we refer the reader to \cite{SV08}.

\begin{figure}[!ht]
\centering
\includegraphics[scale=0.42,clip]{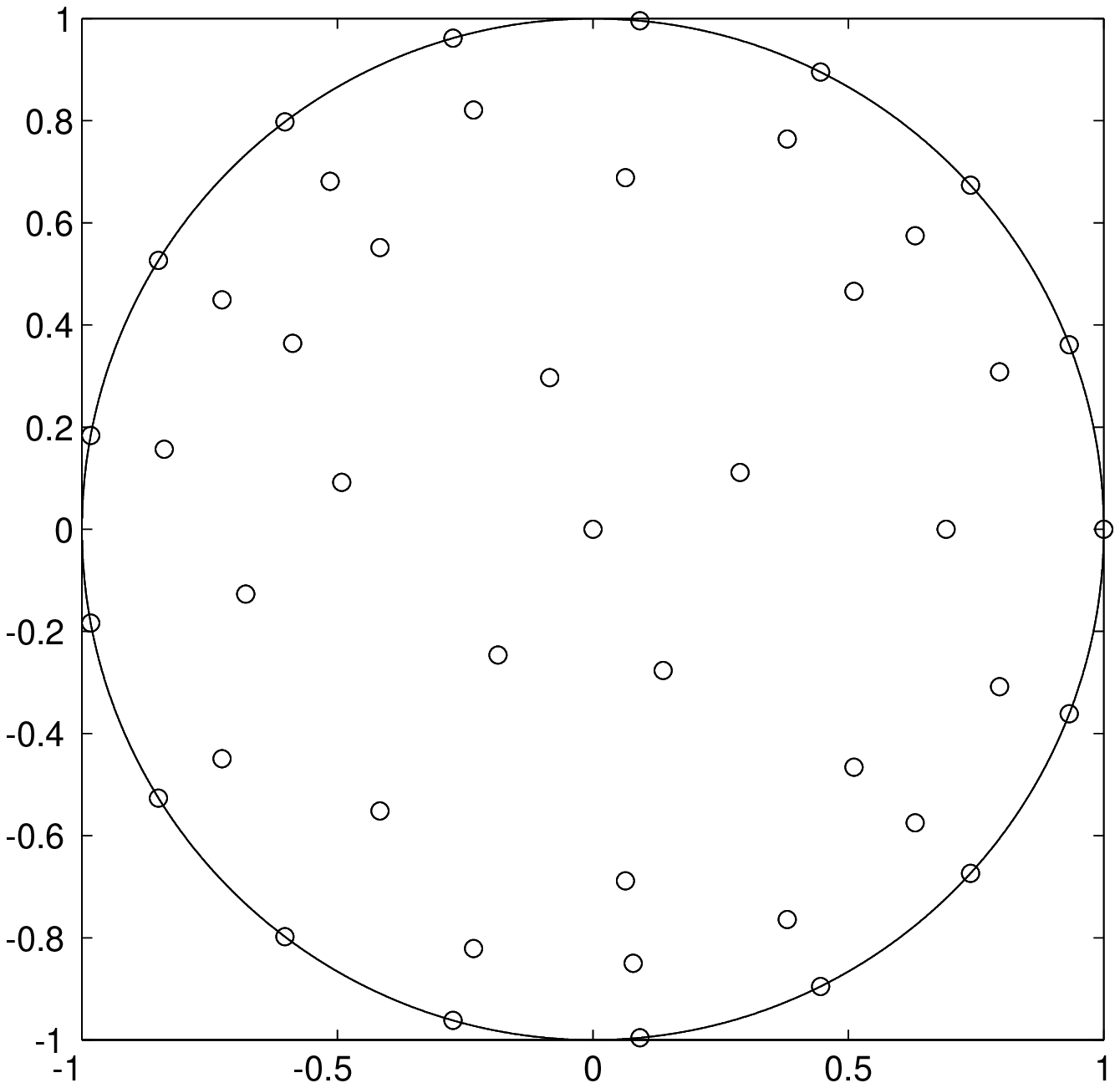}\hfill
\includegraphics[scale=0.42,clip]{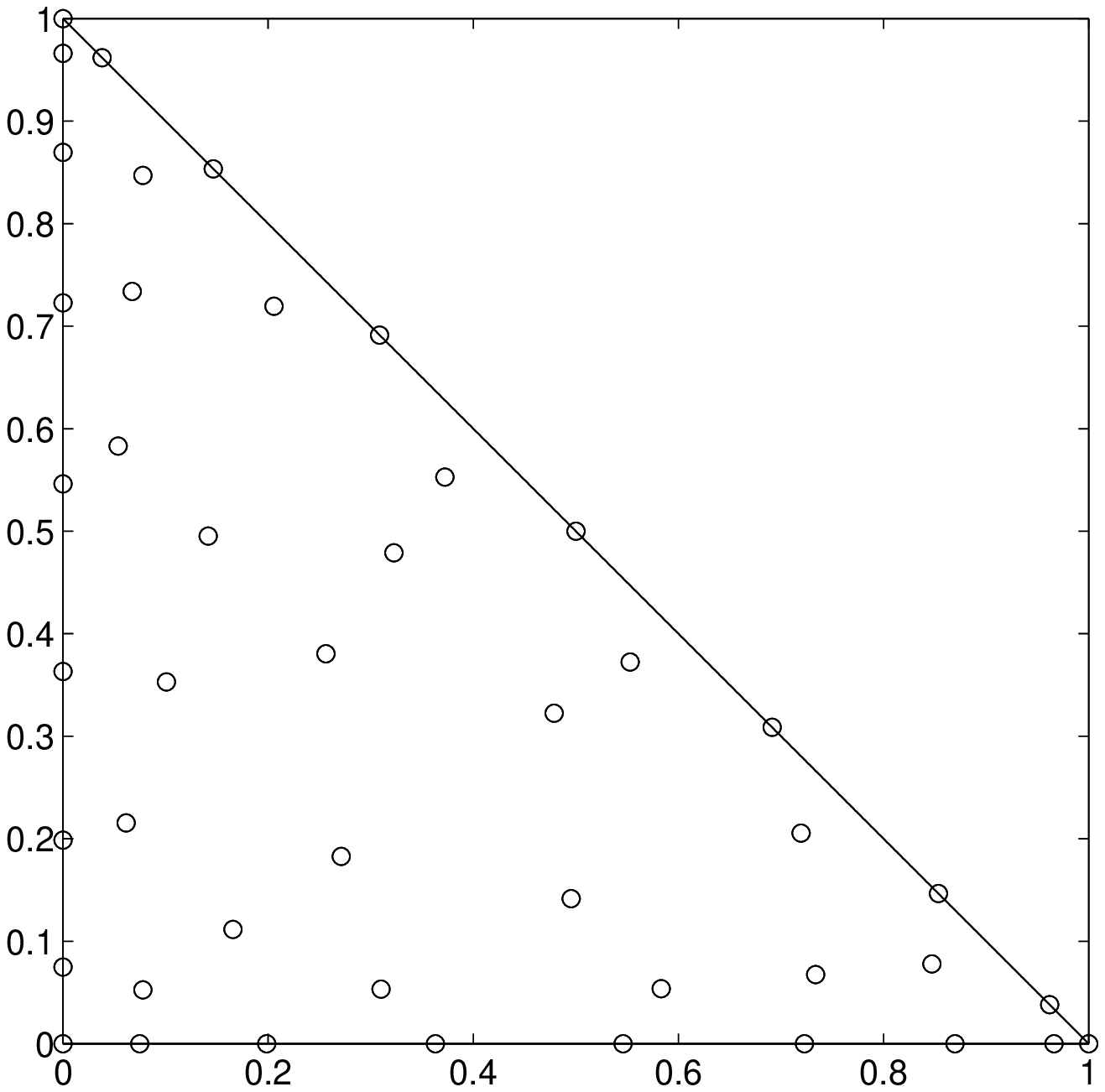}\hfill
\caption{The $45$ Approximate Fekete Points of degree $n=8$
extracted from the geometric WAMs for the disk and the simplex.}
\label{fekDS}
\end{figure}

\begin{figure}[!ht]
\centering
\includegraphics[scale=0.42,clip]{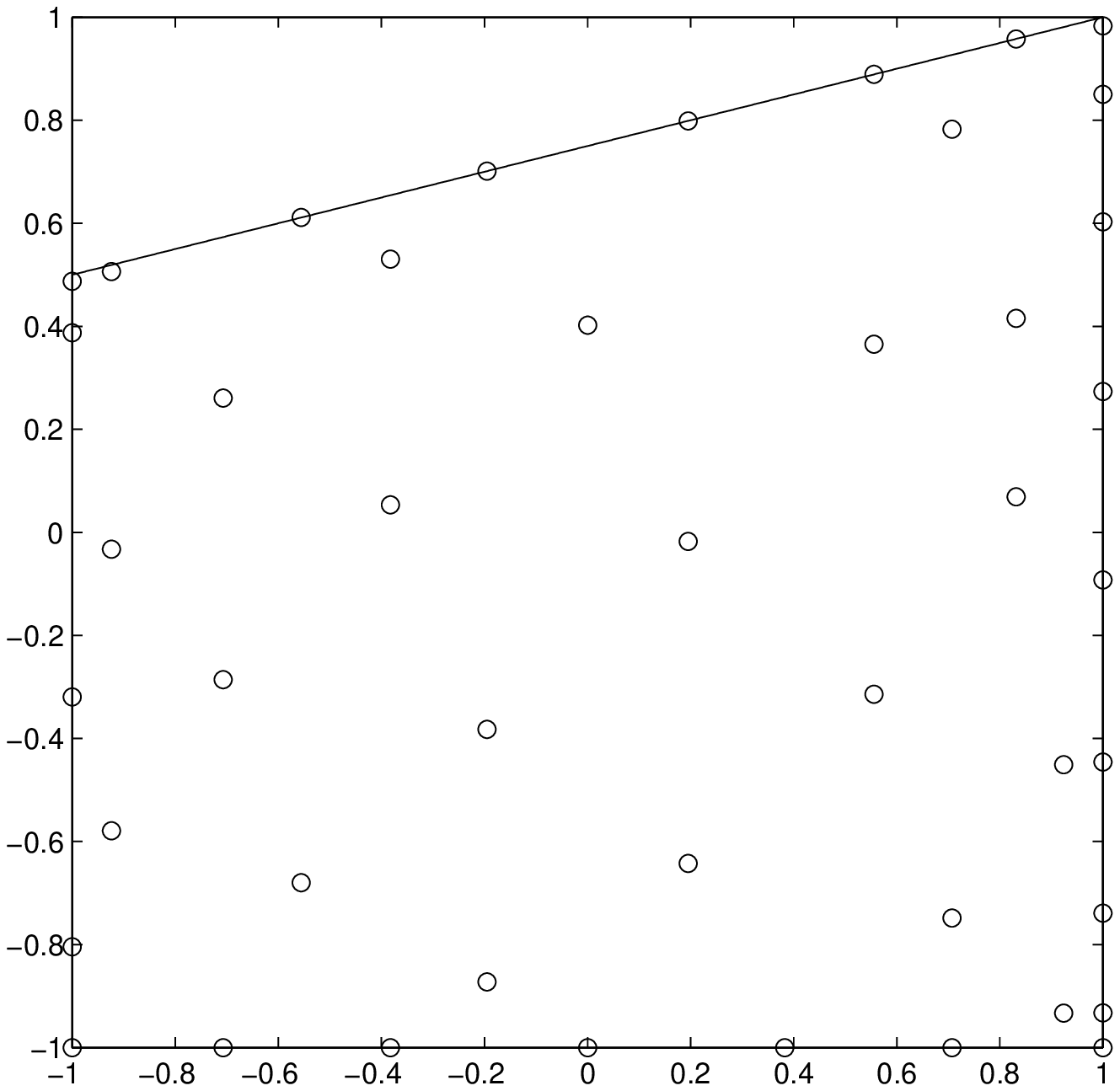}\hfill
\includegraphics[scale=0.42,clip]{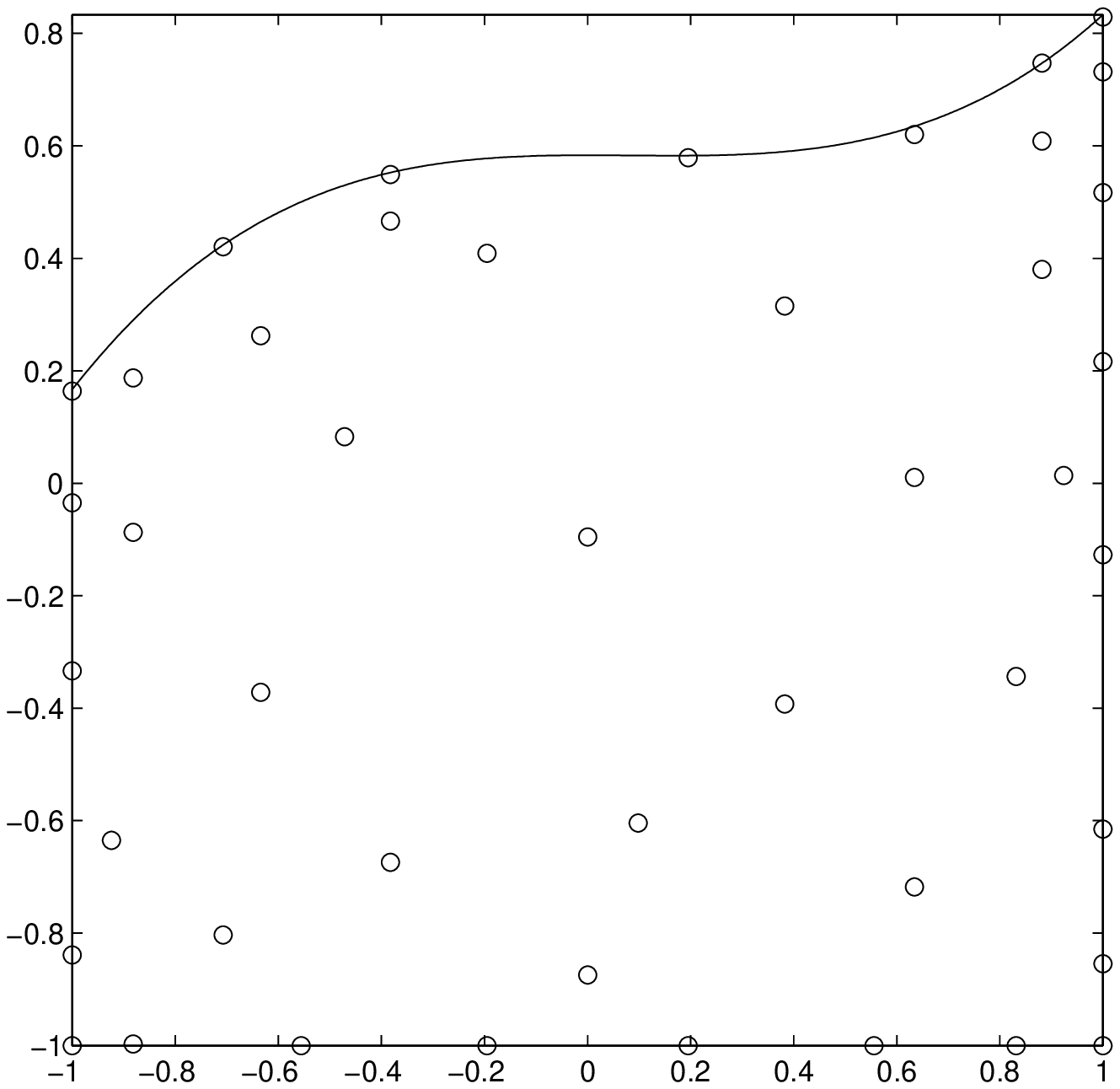}\hfill   
\caption{The $45$ Approximate Fekete Points of degree $n=8$
extracted from the geometric WAMs for a linear and a cubic trapezoid.}
\label{fektrap}
\end{figure}

\begin{figure}[!ht]
\centering
\includegraphics[scale=0.42,clip]{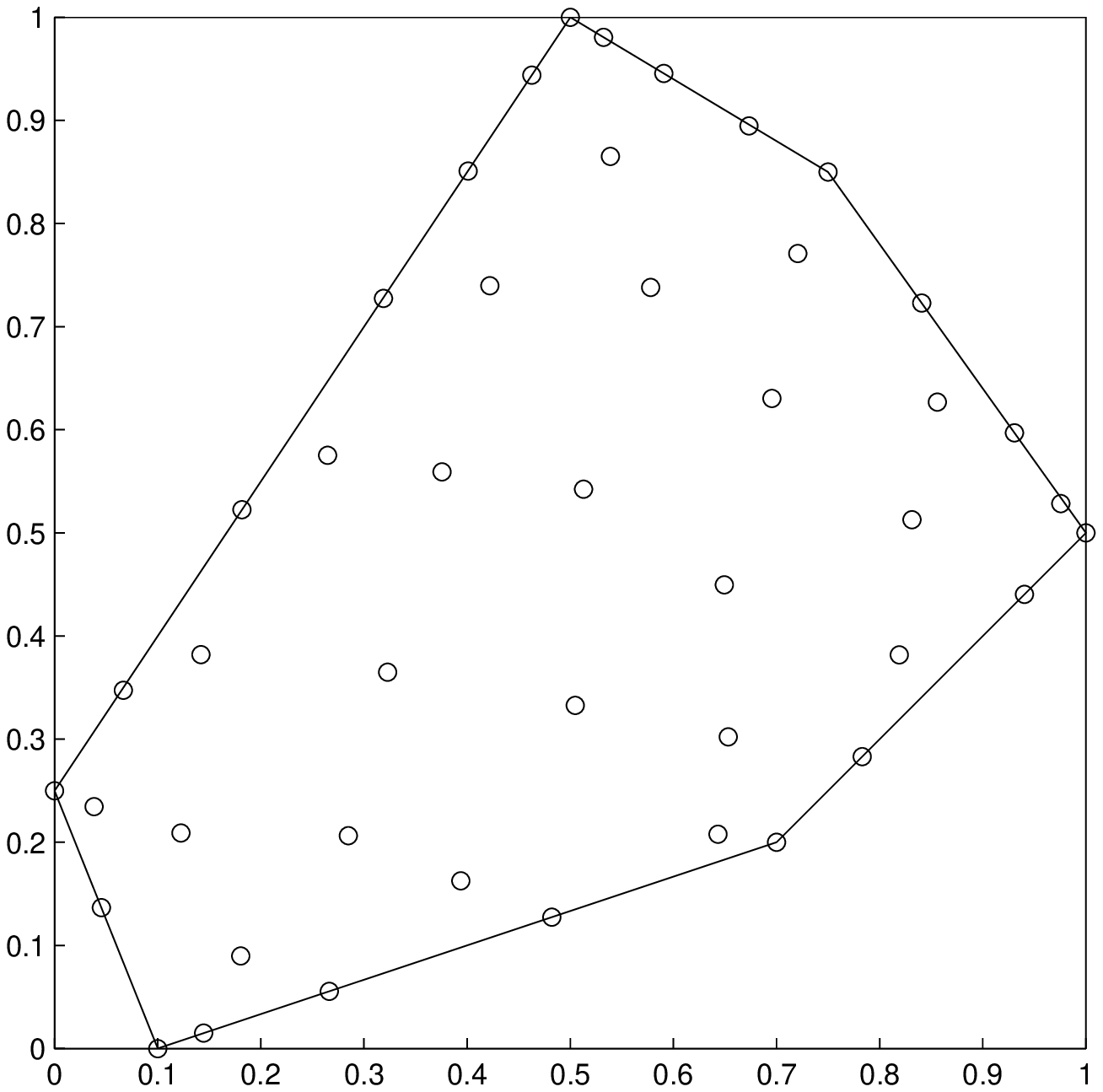}\hfill
\includegraphics[scale=0.42,clip]{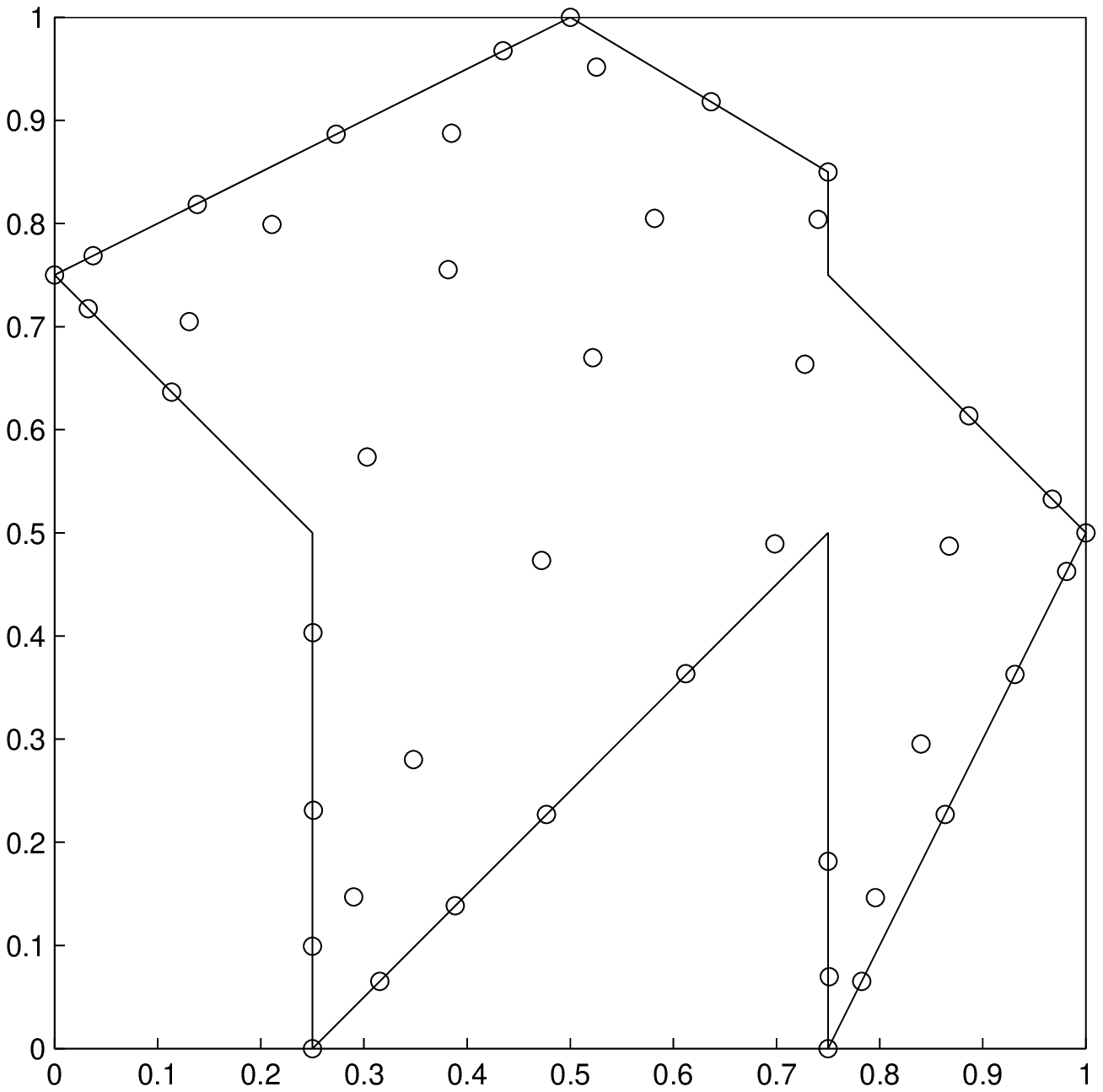}\hfill
\caption{The $45$ Approximate Fekete Points of degree $n=8$
extracted from the geometric WAMs for a convex and
a nonconvex polygon.}
\label{fekP}
\end{figure}

In Figures \ref{fekDS}-\ref{fekP} we display the Approximate Fekete Points 
of degree $n=8$ extracted form the geometric WAMs of Figures 
\ref{disk}-\ref{polygons}. The computational advantage of working with 
a Weakly Admissible Mesh instead of an Admissible Mesh is shown in Table 
1, where we show the cardinalities of the relevant discrete sets in the 
case of the disk. We recall that, following \cite[Thm.5]{CL08}, it is 
always 
possible to construct an AM for a real $d$-dimensional 
compact set that admits a Markov polynomial inequality with exponent $r$, 
by intersection with a uniform grid of stepsize ${\cal O}(n^{-rd})$. 
In particular, convex compact sets have $r=2$, and it is easily seen  
from the Markov polynomial inequality $\|\nabla p(x)\|_2 \leq n^2 
\|p\|_\infty$ valid for every $p\in \bb{P}_n^2$ and $x$ in the disk, 
and the proof of 
\cite[Thm.5]{CL08}, 
that it is sufficient to take a stepsize
$h<1/n^2$. The cardinality of the corresponding AM is then ${\cal 
O}(n^4)$, to be compared with the ${\cal O}(n^2)$ cardinality of 
the geometric WAM (see Section 2.1). For example, at degree $n=30$ an AM 
for the disk 
has more than 2 millions points, whereas the geometric WAM has less than 2 
thousands points. 

\begin{table}[!ht]
\caption{Cardinalities of different point sets in the unit disk (AFP =
Approximate Fekete Points).}
\begin{center}  
\vskip0.1cm   
\begin{tabular}{c c c c c c c}
points & $n=5$ & $n=10$ & $n=15$ & $n=20$ & $n=25$ & $n=30$\\
\hline
AM & 2032 & 31700 & 159692 & 503868 & 1229072 & 2547452\\
WAM & 60 & 220 & 480 & 840 & 1300 & 1860\\
AFP & 21 & 66 & 136 & 231 & 351  & 496\\
\hline
\end {tabular}
\end{center}
\end{table}

In Table 2, we show the Lebesgue constants of the Approximate 
Fekete Points extracted from the geometric WAMs for the disk 
at a sequence of degrees, using the algorithm above with different 
polynomial bases. Such Lebesgue constants have been evaluated 
numerically on a suitable control mesh, much finer than the 
extraction mesh. Without refinement iterations, the best results 
are obtained with the Logan-Shepp basis, which, as is well known, is orthogonal for 
standard Lebesgue measure (cf. \cite{LS75}). On the contrary, with the 
monomial 
basis we face severe 
ill-conditioning and even numerical rank 
deficiency of the Vandermonde matrix, and we get the worst Lebesgue 
constants. After two refinement iterations, 
however, we are working in practice with a discrete orthogonal basis, 
the corresponding Vandermonde matrices are not ill-conditioned,  
and the Lebesgue constants improve and stabilize. Observe that their 
growth is much slower than that of the theoretical bound 
(\ref{leb_fek}). 

\begin{table}[!ht]
\caption{Numerically evaluated Lebesgue constants (nearest integer)
of the Approximate
Fekete Points extracted from the geometric WAM in the unit disk,
with different bases
(Mon = monomial, Che = product Chebyshev, LoS = Logan-Shepp; in
parentheses the number of refinement iterations); $\ast$ means that
the rectangular Vandermonde matrix (\ref{vander}) is numerically
rank-deficient.}
\begin{center}
\vskip0.1cm
\begin{tabular}{c c c c c c c}
basis & $n=5$ & $n=10$ & $n=15$ & $n=20$ & $n=25$ & $n=30$\\
\hline
Mon(0) & 7 & 21 & 34 & 869 & $\ast$ & $\ast$ \\
Mon(2) & 5 & 24 & 32 & 42 & 60 & 81\\
Che(0) & 9 & 23 & 30 & 91 & 1321 & $\ast$\\
Che(2) & 5 & 24 & 32 & 42 & 60 & 81\\
LoS(0) & 7 & 20 & 32 & 52 & 87 & 119\\
LoS(2) & 5 & 24 & 32 & 42 & 60 & 81\\
\hline
\end {tabular}
\end{center}
\end{table}

In Table 3 we compare, for the three test functions below, 
the errors (in the uniform norm) of discrete least 
squares approximation on the AM and on the WAM of the disk, and  
of interpolation at the Approximate Fekete Points extracted 
from the WAM (with two refinement 
iterations). 
The test functions exhibit different regularity: the first is 
analytic entire, the second is analytic nonentire (a bivariate version
of the classical Runge function), the third is $C^1$ but has a singularity 
of the second derivatives at the origin. 

\begin{itemize}
\item {\bf test function 1:} $f(x_1,x_2)=\cos{(x_1+x_2)}$
\item {\bf test function 2:} $f(x_1,x_2)=1/(1+16(x_1^2+x_2^2))$
\item {\bf test function 3:} $f(x_1,x_2)=(x_1^2+x_2^2)^{3/2}$
\end{itemize}
\noindent
The Vandermonde 
matrices have been constructed using the product Chebyshev basis.  
Both the least squares and the interpolation polynomial coefficients 
have been computed by the standard Matlab ``backslash'' solver, and 
the errors have been evaluated on a suitable control mesh. 

\begin{table}[!ht]
\caption{Uniform errors of polynomial approximations on different point
sets in the unit disk, for the three test functions above; $\ast$
means computational failure due to large dimension (see Table 1).}
\begin{center}
\vskip0.1cm
\begin{tabular}{c c c c c c c c}
 & points & $n=5$ & $n=10$ & $n=15$ & $n=20$ & $n=25$ & $n=30$\\
\hline
test 1 & LS AM & 9E-4 & 3E-10 & $\ast$ & $\ast$ & $\ast$ & $\ast$\\
 & LS WAM & 5E-4 & 1E-10 & 3E-15 & 7E-15 & 6E-15 & 2E-14\\
 & interp AFP & 1E-3 & 3E-10 & 2E-15 & 2E-15 & 2E-15 & 3E-15\\
\hline
test 2 & LS AM & 4E-1 & 1E-1 & $\ast$ & $\ast$ & $\ast$ & $\ast$\\   
 & LS WAM & 5E-1 & 7E-2 & 5E-2 & 6E-3 & 4E-3 & 5E-4\\
 & interp AFP & 5E-1 & 7E-2 & 5E-2 & 6E-3 & 4E-3 & 5E-4\\
\hline

test 3 & LS AM & 2E-2 & 2E-3 & $\ast$ & $\ast$ & $\ast$ & $\ast$\\
 & LS WAM & 2E-2 & 1E-3 & 7E-4 & 1E-4 & 2E-4 & 4E-5\\
 & interp AFP & 2E-2 & 1E-3 & 7E-4 & 1E-4 & 2E-4 & 4E-5\\
\hline
\end {tabular}
\end{center}
\end{table}

Notice that with the AM we have computational failure (``out of 
memory'') in our computing system already at degree $n=15$, due to the large 
cardinality of the discrete set; see Table 1. The least squares error  
on the WAM is close to that on the AM (when comparable), which shows 
that geometric WAMs are a good choice for polynomial approximation, 
with a low computational cost. It is worth observing that in the theoretical 
estimate (\ref{LS}) we even have  $C(A_n)\sqrt{{\rm 
card}(A_n)}={\cal O}(n^2)$ for the AM, and $C(A_n)\sqrt{{\rm 
card}(A_n)}={\cal O}(n \log^2{(n)})$ for the WAM, but we recall 
that these are overestimates, the term $\sqrt{{\rm card}(A_n)}$ being 
in some way ``artificial'' (cf. \cite[Thm.2]{CL08}). 

The following tables are devoted to numerical tests on the 
other domains. First, in Table 4 we show the Lebesgue constants of the 
Approximate Fekete Points extracted from the geometric WAMs described in 
Section 2. Again, the  growth is much slower than that of the theoretical 
bound (\ref{leb_fek}).

As for the 
simplex, the Lebesgue constant of 
our Approximate Fekete Points points is larger than 
that of the best points known in the literature. The case of the simplex 
has been widely studied and several specialized approaches have been 
proposed for the 
computation of Fekete or other good interpolation points, due to the 
relevance in the numerical treatment of PDEs by spectral-element and 
high order finite-element methods: see e.g. \cite{HW08,PR06,TWV00,W06} and
references therein. 

On the other hand, our method for computing 
Approximate Fekete Points via geometric WAMs is quite general and 
flexible, since it allows to work on a wide class of compact sets and, 
differently from other computational approaches, up to  
reasonably high interpolation degrees.       
The good quality of the discrete sets used for least squares 
approximation and for interpolation is evidenced by Tables 5-9. 
We observe that the singularity of the third 
test function is in the interior of the domains, apart from the 
simplex where it is located at a vertex (where the discrete points 
cluster). This explains the better results with the simplex for 
this function. In the case of the two polygons, a change of variables 
is made 
in order to put the problem in the reference square $[-1,1]^2$.   

The availability of good interpolation points in compact sets with 
various geometries has a number of potential applications. One for 
example is 
connected to numerical cubature. Indeed, when the moments of the 
underlying polynomial basis are known \cite{SV07,SV09}, cubature weights 
associated to 
the Approximate Fekete Points can 
be computed 
as a by-product of the algorithm, simply by using the moments 
vector as right-hand 
side $b$. This gives an algebraic cubature formula, that can 
be used directly, or as a starting point towards the computation 
of a minimal formula, by the method developed in \cite{TWB07}. 

Another relevant application concerns 
the numerical treatment of PDEs, where a renewed interest is arising 
in global polynomial methods, such as collocation and discrete least 
squares methods, over general domains (see, e.g., \cite{MNPP09}). 
Recently, Approximate Fekete 
Points have been successfully used for discrete least squares 
discretization of elliptic equations \cite{Z09}. Moreover, again in the 
context of numerical PDEs, Approximate Fekete Points for polygons could play 
a role in connection with discretization methods on polygonal (non 
simplicial) meshes (see, e.g., \cite{SM06} and references therein). 

\begin{table}[!ht]
\caption{Numerically evaluated Lebesgue constants (nearest integer) 
of the Approximate 
Fekete Points extracted from the geometric WAMs, on different compact 
sets (product Chebyshev basis with 2 refinement iterations).}     
\begin{center}                                                   
\vskip0.1cm                                        
\begin{tabular}{c c c c c c c}
set & $n=5$ & $n=10$ & $n=15$ & $n=20$ & $n=25$ & $n=30$\\
\hline
disk & 5 & 24 & 32 & 42 & 60 & 81\\
simplex & 5 & 15 & 25 & 48 & 62 & 80\\
linear trap & 6 & 19 & 34 & 37 & 31 & 54\\
cubic trap & 6 & 16 & 35 & 45 & 40 & 75\\
conv polyg  & 7 & 13 & 22 & 53 & 53 & 66\\
nonconv polyg & 5 & 18 & 36 & 35 & 45 & 80\\
\hline
\end {tabular}
\end{center}  
\end{table}   

\begin{table}[!ht]
\caption{Uniform errors of polynomial approximations in the unit simplex, 
for the three test functions above.}
\begin{center}
\vskip0.1cm
\begin{tabular}{c c c c c c c c}
 & points & $n=5$ & $n=10$ & $n=15$ & $n=20$ & $n=25$ & $n=30$\\
\hline
test 1 
 & LS WAM & 7E-7 & 8E-15 & 3E-15 & 4E-15 & 4E-15 & 6E-15\\
 & interp AFP & 2E-6 & 2E-14 & 1E-15 & 3E-15 & 3E-15 & 5E-15\\
\hline
test 2 
 & LS WAM & 2E-2 & 5E-4 & 1E-5 & 4E-7 & 1E-8 & 4E-10\\
 & interp AFP & 5E-2 & 2E-3 & 4E-5 & 2E-6 & 3E-8 & 2E-9\\
\hline
test 3 
 & LS WAM & 7E-4 & 5E-6 & 4E-7 & 8E-8 & 2E-8 & 7E-9\\
 & interp AFP & 8E-4 & 2E-5 & 1E-6 & 2E-7 & 6E-8 & 3E-8\\
\hline
\end {tabular}
\end{center}
\end{table}

\begin{table}[!ht]
\caption{As in Table 5, for the linear 
trapezoid of Figure 3.}
\begin{center}
\vskip0.1cm
\begin{tabular}{c c c c c c c c}
 & points & $n=5$ & $n=10$ & $n=15$ & $n=20$ & $n=25$ & $n=30$\\
\hline
test 1
 & LS WAM & 3E-3 & 5E-9 & 1E-13 & 3E-15 & 4E-15 & 9E-15\\
 & interp AFP & 8E-3 & 2E-8 & 3E-13 & 4E-15 & 3E-15 & 4E-15\\
\hline
test 2
 & LS WAM & 2E-1 & 2E-1 & 1E-1 & 3E-2 & 1E-2 & 5E-3\\
 & interp AFP & 3E-1 & 2E-1 & 2E-1 & 3E-2 & 2E-1 & 1E-2\\
\hline
test 3
 & LS WAM & 3E-2 & 4E-3 & 2E-3 & 5E-4 & 2E-4 & 1E-4\\
 & interp AFP & 5E-2 & 4E-3 & 3E-3 & 5E-4 & 3E-4 & 2E-4\\
\hline
\end {tabular}
\end{center}
\end{table}

\begin{table}[!ht]
\caption{As in Table 5, for the cubic
trapezoid of Figure 3.}
\begin{center}                      
\vskip0.1cm   
\begin{tabular}{c c c c c c c c}
 & points & $n=5$ & $n=10$ & $n=15$ & $n=20$ & $n=25$ & $n=30$\\
\hline
test 1
 & LS WAM & 2E-3 & 6E-9 & 6E-14 & 3E-15 & 4E-15 & 5E-15\\
 & interp AFP & 6E-3 & 1E-8 & 1E-13 & 5E-15 & 3E-15 & 4E-15\\
\hline
test 2
 & LS WAM & 4E-1 & 2E-1 & 6E-2 & 3E-2 & 9E-3 & 5E-3\\ 
 & interp AFP & 5E-1 & 2E-1 & 7E-2 & 5E-2 & 1E-2 & 6E-3\\
\hline
test 3
 & LS WAM & 3E-2 & 3E-3 & 9E-4 & 5E-4 & 2E-4 & 2E-4\\
 & interp AFP & 6E-2 & 5E-3 & 9E-4 & 7E-4 & 2E-4 & 2E-4\\
\hline
\end {tabular}
\end{center}  
\end{table}

\begin{table}[!ht]
\caption{As in table 5 for the convex polygon 
of Figure 4 and the three test functions $f(2x_1-1,2x_2-1)$ above.}
\begin{center}
\vskip0.1cm
\begin{tabular}{c c c c c c c c}
 & points & $n=5$ & $n=10$ & $n=15$ & $n=20$ & $n=25$ & $n=30$\\
\hline
test 1 & LS WAM & 7E-4 & 1E-9 & 7E-15 & 9E-15 & 1E-14 & 2E-14\\
 & interp AFP & 1E-3 & 4E-9 & 6E-15 & 4E-15 & 4E-15 & 5E-15\\
\hline
test 2 & LS WAM & 4E-1 & 1E-1 & 4E-2 & 2E-2 & 4E-3 & 1E-3\\
 & interp AFP & 5E-1 & 1E-1 & 4E-2 & 2E-2 & 9E-3 & 3E-3\\
\hline
test 3 & LS WAM & 2E-2 & 2E-3 & 6E-4 & 3E-4 & 1E-4 & 9E-5\\
 & interp AFP & 2E-2 & 2E-3 & 6E-4 & 3E-4 & 1E-4 & 8E-5\\
\hline
\end {tabular}
\end{center}
\end{table}

\begin{table}[!ht]
\caption{As in Table 8 for the nonconvex 
polygon 
of Figure 4.}
\begin{center}
\vskip0.1cm
\begin{tabular}{c c c c c c c c}
 & points & $n=5$ & $n=10$ & $n=15$ & $n=20$ & $n=25$ & $n=30$\\
\hline
test 1 & LS WAM & 5E-4 & 3E-10 & 1E-14 & 2E-14 & 3E-14 & 4E-13\\
 & interp AFP & 6E-4 & 5E-10 & 3E-15 & 3E-15 & 3E-15 & 4E-15\\
\hline
test 2 & LS WAM & 4E-1 & 2E-1 & 5E-2 & 2E-2 & 5E-3 & 1E-3\\   
 & interp AFP & 6E-1 & 2E-1 & 5E-2 & 2E-2 & 5E-3 & 2E-3\\
\hline
test 3 & LS WAM & 2E-2 & 3E-3 & 7E-4 & 3E-4 & 1E-4 & 9E-5\\
 & interp AFP & 4E-2 & 3E-3 & 8E-4 & 3E-4 & 1E-4 & 7E-5\\
\hline
\end {tabular}
\end{center}
\end{table}

\end{document}